\documentclass{article}
\usepackage{amsmath,amsfonts,amssymb,array,tabularx,enumerate}
\usepackage{epigraph}
\usepackage{enumitem}      
\usepackage{arxiv}
\usepackage{float}
\usepackage[utf8]{inputenc} 
\usepackage{fvextra} 
\usepackage{capt-of}   
\usepackage[T1]{fontenc}    
\usepackage{verbatim}
\DeclareUnicodeCharacter{03B1}{\ensuremath{\alpha}}
\usepackage{hyperref}       
\usepackage{url}  
\usepackage{booktabs}
\usepackage{amsfonts}  
\usepackage{nicefrac}
\usepackage{microtype}  
\usepackage{lipsum} 
\usepackage{graphicx}
\usepackage{doi}
\usepackage[numbers,sort&compress]{natbib}

\usepackage{xcolor}
\definecolor{ghl}{rgb}{1.0, 0.84, 0.0}  

\newtheorem{thm}{Theorem}[section]
\newtheorem{cor}[thm]{Corollary}
\newtheorem{lem}[thm]{Lemma}

\title{Descriptions of Cantor Sets:\\ A Set‑Theoretic Survey and Open Problems}
\author{Mohsen Soltanifar\thanks{Contact: mohsen.soltanifar[at]alumni.utoronto.ca.} \href{https://orcid.org/0000-0002-5989-0082}{\includegraphics[scale=0.06]{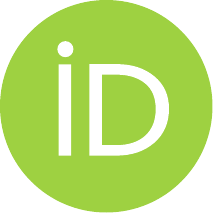}}
\\ Biostatistics Division, Dalla Lana School of Public Health, University of Toronto\\
620-155 College Street, Toronto, ON M5T 3M7, Canada}
 
\renewcommand{\shorttitle}{\textit{arXiv} Template}

\begin{document}
\maketitle

\begin{abstract}
This survey synthesizes the principal descriptive set-theoretic perspectives on deterministic Cantor sets on the real line and charts directions for future study. After recounting their historical genesis and compiling an up-to-date taxonomy, we review the Borel hierarchy and four hierarchically ordered representations—general, nested, iterated-function-system (IFS), and $q-$ary expansion—presented from the most general to the most specific set-theoretic description of deterministic Cantor sets. We then present explicit closed-form descriptions for two thin families of measure-zero Cantor sets and, for the augmented ``thick'' family of positive Lebesgue measure, a new closed recursive formula for the endpoint maps governing every retained interval at every finite stage; we further show that the classical middle-third set lies in the intersection of all three families, providing a unifying example across the four representations.The survey closes by isolating several open problems in four directions, aiming to provide mathematicians with a coherent platform for further descriptive set-theoretic investigations into Cantor-type sets on the real line.
\end{abstract}

\vspace{1em}

\epigraph{``A set is a Many that allows itself to be thought of as a One.''}{--- Georg Cantor (1845--1918)}

\section{Introduction}\label{1.Int}
The Introduction is structured as follows. Section \ref{1.1.HB} reviews the historical development of Cantor sets, while Section \ref{1.2.TCS} presents an up-to-date taxonomy of their principal variants. Section \ref{1.3.MS} sets out the motivation for undertaking a set-theoretic survey of these constructions, and Section \ref{1.4.OP} closes the Introduction with a synopsis of the paper’s overall organization.\par 

\subsection{Historical Background}\label{1.1.HB}
The conceptual roots of the Cantor Sets lie in  Georg Cantor's 1872 introduction of the \emph{derived set} (the set of all limit points of a given subset of $\mathbb{R}$), which laid the groundwork for perfect, nowhere-dense subsets of $\mathbb{R}.$ Shortly thereafter, H.~J.~S.~Smith’s 1875 note provided the first printed construction of a
“Cantor‐type’’ set by iteratively excising subintervals of \([0,1]\), thereby demonstrating that a closed set can be uncountable yet possess arbitrarily small outer content, although the significance of his example went largely unnoticed at the time.\, Cantor returned to the theme in his six–part memoir \emph{Über unendliche lineare Punktmannichfaltigkeiten} (1879–1884); Part V (written October 1882, published 1883) contains the now‐classical \emph{3-ary or ternary expansion set} (later called the middle-third Cantor set), presented as a perfect set that is nowhere dense.\, In a letter dated November 1883 Cantor extends this construction of the Cantor set by defining the \emph{Cantor function}, a continuous, non‐decreasing map that is constant on the removed intervals and
thus furnishes a counter-example to naïve forms of the Fundamental Theorem of Calculus.\, His 1884 paper \emph{De la puissance des ensembles parfaits de points} in \emph{Acta Mathematica} analyzed the cardinalities of perfect sets, cementing the role of the ternary Cantor set as both a set-theoretic and analytic paradigm. These milestones—Smith’s overlooked precursor, Cantor’s derived-set framework, the explicit ternary construction, and the analytic Cantor function—together mark the emergence of what is now simply known as the Cantor set \cite{Fleron1994}.\par

\subsection{Taxonomy of Cantor Sets}\label{1.2.TCS}
Various classes and approaches to the construction of Cantor sets have been extensively studied, each generating distinct subclasses and mathematical structures. As of the date of this survey (31 January 2026), mathematicians have identified numerous specific types of Cantor sets, including Affine, Autonomous, Balanced, Bunched, Central, Cookie-Cutter, Deterministic, Discrete, Distorted, Dusts, Fat, Generalized, Genus-g, Geometric, Hairy, High Dimensional, Homogeneous, Hyperbolic, Iteration Function System (IFS) type, Invariant, Kinetic, Linear, Minimal, Multi-model, Nearly-affine, Non autonomous, Non-homogeneous, Non-archimedean (p-adic), Non-central, Non-Symmetric, Nonstandard, p-Adic Heisenberg, Random, Regular, Scrawny, Self-Similar, Semi-bounded, Smooth, Sticky, Superior, Smith-Volterra-Cantor(SVC) type, Symbolic, Symmetric, Tame, Thick, Thin, Twofold, Ultrametric, Uniform, Universal, Visible, and Wild. Comprehensive references for these subclasses can be found in \cite{Edgar2007,Valline2013,Falconer2014,Gorodetski2015,Lin2024,
McLoughlin2014,Jumah2024,Nassiri2025,Nowakowski2021,Batakis2015,
Batakis2024,Bellissard2013,Hieronymi2018,Hu2008,Hassan2018,Garcia2011,
Fillman2023,Denette2016,Tresser1991,Moreira2020,Lapidus2008,
Cheraghi2019,Cabrelli2011,Cockerill2002,Hare2015,Zeng2017,
Zhou2019,Honary2006,Fletcher2014,Akhadkulov2019,Babich1992,
Krushkal2018,Kamalutdinov2018,Kong2014,Semmes2011,Li2008,
Martens1994,Nakata1997,Rani2010,Fletcher2024,Eswarathasan2023,
Hare2024,Bartoszewicz2017,Mantica2015}. In this paper, however, our primary focus is on deterministic Cantor sets defined explicitly on the real line, and we direct the interested reader to relevant literature for more general or abstract variants.\par 

\subsection{Motivation and Scope}\label{1.3.MS}
Since Georg Cantor’s original construction in the 1880s, the Cantor set has become a touchstone throughout mathematics and beyond. On the real line it provides one of the simplest explicit exemplars of a fractal: it is  nowhere dense (its closure has empty interior), totally disconnected
(its only connected subsets are singletons), perfect (closed with no isolated points), and uncountable, yet its Lebesgue measure can be tuned from zero (the classical middle-third set) to any prescribed value in $(0,1).$ Such extremal combinations of topological and measure-theoretic properties make Cantor sets a natural testing ground for  ``pathological" phenomena—most famously the Cantor–Lebesgue function, a continuous, non-decreasing map that is constant on a set of full measure. In constructive and descriptive settings the Cantor set reappears as a paradigmatic building block: for example, basic Cantor-type function blocks span 5/28 of the canonical representatives in the function space $F(R,R),$ \cite{Soltanifar2023} and they appear to prove the existence of aleph-two deterministic as well as aleph-two random fractals on the real line \cite{Soltanifar2021,Soltanifar2022}. Outside pure mathematics, Cantor-type patterns have been cited in art—early Egyptian column motifs \cite{Frame2025,Kainth2023}, in modeling financial markets \cite{Jin2021} —and in nature, such as ring-like stratification of Saturn when idealized as repeated circular products \cite{Mandelbrot1982}.\par

\paragraph{Applications and Cross-Disciplinary Connections} Beyond the illustrations just listed, Cantor-type sets recur across mathematics as both structural prototypes and technical tools, and a brief tour of these connections clarifies why a unified descriptive survey is timely. We summarize seven of their principal connections, several of which underlie the open problems posed in Section~\ref{4.D}.

\begin{enumerate}

\item[ACDC1.] \textbf{Dynamical systems and symbolic dynamics.}
Hyperbolic invariant sets of smooth dynamical systems frequently carry the topology of a Cantor set, and the conjugacy between such systems and full shifts on finite alphabets is the foundational technique of symbolic dynamics~\cite{Lind1995}. The horseshoe construction of Smale and the wild hyperbolic sets of Newhouse~\cite{Newhouse1979} are paradigmatic examples, and the IFS attractors discussed in Section~\ref{2.STF} are their simplest one-dimensional analogues.

\item[ACDC2.] \textbf{Number theory and Diophantine approximation.}
The set of badly approximable real numbers in $[0,1]$ is, up to measure zero, a countable union of Cantor-type sets indexed by partial-quotient bounds. The arithmetic theory of sums and differences of Cantor sets---initiated by Newhouse and developed extensively by Palis, Takens, and Moreira--Yoccoz~\cite{Palis1993, Moreira2001}---asks when $\Gamma + \Gamma'$ or $\Gamma - \Gamma'$ contains an interval, with applications to homoclinic bifurcations and to the Markov and Lagrange spectra~\cite{Moreira2020}.

\item[ACDC3.] \textbf{Harmonic analysis.}
Self-similar measures supported on Cantor sets serve as test cases and counterexamples in Fourier analysis. Salem sets---compact sets whose Hausdorff and Fourier dimensions coincide---are constructed using randomized Cantor-type schemes~\cite{Salem1951}, and Cantor measures appear centrally in modern work on arithmetic progressions in sparse sets~\cite{Laba2009}.

\item[ACDC4.] \textbf{Geometric measure theory and fractal geometry.}
Cantor sets are the canonical examples in the theory of Hausdorff dimension, self-similar measures, and rectifiability~\cite{Mattila1995, Falconer2014}. The Hausdorff-Dimension Theorem (cf.\ Section~\ref{3.DFCS}) is itself most transparently proved by exhibiting an explicit one-parameter family of Cantor sets, and the recent extensions in~\cite{Soltanifar2021, Soltanifar2022} place such constructions in a broader cardinality-theoretic framework.

\item[ACDC5.] \textbf{Descriptive set theory and topology.}
Brouwer's classical theorem characterizes the Cantor space $\{0,1\}^{\mathbb{N}}$ as the unique (up to homeomorphism) non-empty perfect, compact, totally disconnected, metrizable space~\cite{Brouwer1910, Kechris1995}. Every Cantor set on the real line is therefore homeomorphic to every other, which makes the descriptive distinctions surveyed in this paper---arising not from topology but from Lebesgue measure, Hausdorff dimension, and explicit arithmetization---the more striking.

\item[ACDC6.] \textbf{Probability theory.}
Random analogues of the constructions in Section~\ref{2.STF} give rise to random recursive fractals~\cite{Mauldin1986, Soltanifar2022} and to determinantal point processes whose intensity measures are concentrated on generalized Cantor sets~\cite{Lin2024}. These constructions are the stochastic counterpart of the deterministic survey pursued here.

\item[ACDC7.] \textbf{Further applications outside pure mathematics.}
In addition to the art-, finance-, and physics-related references cited above, Cantor-type structures have been invoked as descriptive models in signal processing and information theory ~\cite{Schroeder2009}.\par  

Several of these directions---particularly the arithmetic of intersections in dynamical systems and the harmonic-analytic study of Cantor measures---motivate the open problem (OP5) on common intersections posed in Section~\ref{4.D}.\par
\end{enumerate}

Despite this wide relevance, no unified survey has yet cataloged the diverse set-theoretic descriptions of Cantor sets, nor identified which directions are fully developed and which remain open. The present paper closes that gap by reviewing known explicit and implicit set-theoretic formulaes on the real line for the deterministic cases, classifying them by structural features, and highlighting avenues for further investigation.  While the book-length treatment of Vallin~\cite{Valline2013} catalogs many properties and applications of Cantor sets, it does not organize their descriptions into a hierarchical set-theoretic framework, nor does it supply explicit endpoint recursions for positive-measure middle-$\alpha$ families;
the present work is, to our knowledge, the first survey to do so.\par 
\paragraph{Novelty and contributions.} Beyond consolidating known material, this survey advances the descriptive set-theoretic study of Cantor sets in four concrete directions viewed here:
\begin{enumerate}
\item[NC1.] \emph{A unified four-tier hierarchy.} We organize the
general, nested, iterated-function-system, and $q$--ary descriptions under a single set-theoretic language, making translation between representations routine.
\item[NC2.] \emph{An explicit endpoint recursion for positive measure.}Theorem~\ref{thm3} supplies a closed recursive formula for the endpoint maps $a_{n,k}(p/q),\, b_{n,k}(p/q)$ of the augmented thick family $\Gamma_3(p/q,2)$ with $0 < p/q \leq 1/3$; to the best of our knowledge this is the first such explicit recursion in the literature.
\item[NC3.] \emph{The middle-third set as triple intersection.}
Corollary~\ref{cormain} identifies the classical Cantor set as
a member of all three families simultaneously, with equivalence across  all four representations made explicit.
\item[NC4.] \emph{Algorithmic reproducibility and an open agenda.} The accompanying \texttt{R} implementation (Appendix~B) realizes the recursion in (N2) for exact finite-stage computation, and Section~\ref{4.D} distills five open problems (OP1--OP5) that the present framework renders amenable to computational experiment.
\end{enumerate}
These contributions distinguish the present work from prior catalogs and position it as a platform for further descriptive set-theoretic research on Cantor-type sets of the real line. \par
\subsection{Outline of the Paper}\label{1.4.OP}
This survey focuses exclusively on deterministic Cantor sets on the real line, examined through a set-theoretic lens, and highlights several descriptive perspectives that remain undiscovered. Section~\ref{2.STF} reviews the necessary background in descriptive set theory: the Borel hierarchy is recalled, and four nested levels of set-theoretic representation—ranging from the broadest class of Cantor sets to the narrowest class of construction—are delineated. Section~\ref{3.DFCS} then supplies descriptive formulas for three  families of Cantor sets; the middle-third Cantor set lies in the intersection of all three, thereby serving as a unifying example. Section~\ref{4.D} enumerates current challenges, states open problems, and sketches directions for future work. An appendix collects intermediate background, and interim results, and provides annotated computer code that automates some of the constructions, offering a platform for subsequent mathematical investigations.\par 

\section{Set‑Theoretic Foundations}\label{2.STF}
This section deals with the set theoretic foundations for description of Cantor sets. We start with fundamental definitions of $F_{\sigma},$ $G_{\delta}$ and the Borel Hierarchy \cite{Moschovakis2009,Kechris1995}. Then, we present the proposed set theoretic  representations of Cantor sets in hierarchy from macro-level to micro-level in four levels: (i) The general representation, (ii) the nested representation, (iii) the Iterated Function System representation, and, (iv) the $q-$ ary expansion. Recurring specialized terms are collected in Appendix~C for the convenience of non-specialists.\par 

\subsection{The $F_{\sigma}$ and $G_{\delta}$ Properties and the Borel Hierarchy} \label{2.1.FGBH}

Let \(\mathcal{B}([0,1])\) denote the Borel \(\sigma\)-algebra on the closed unit interval.  
Set
\[
  \Sigma^{0}_{1}:=\{\,U\subseteq[0,1] \mid U \text{ is open}\},
  \qquad
  \Pi^{0}_{1}:=\{\,F\subseteq[0,1] \mid F^c\in\Sigma^{0}_{1}\}
  \;(=\text{ all closed sets}),
\]
and define the \emph{Borel hierarchy} inductively for \(n\ge 1\) by:
\[
  \Sigma^{0}_{n+1}
  :=\Bigl\{\bigcup_{k=1}^{\infty}A_k \,\Bigm|\, A_k\in\Pi^{0}_{n}\Bigr\},
  \qquad
  \Pi^{0}_{n+1}
  :=\Bigl\{\bigcap_{k=1}^{\infty}B_k \,\Bigm|\, B_k\in\Sigma^{0}_{n}\Bigr\}.
\]
The classes \(\Sigma^{0}_{2}\) and \(\Pi^{0}_{2}\) are traditionally denoted:  
\[
  F_{\sigma}:=\Sigma^{0}_{2}
  \quad (\text{\emph{countable unions of closed sets}}),\qquad
  G_{\delta}:=\Pi^{0}_{2}
  \quad (\text{\emph{countable intersections of open sets}}).
\]
Thus, every closed subset of \([0,1]\) is \(F_{\sigma}\), and every open subset of \([0,1]\) is \(G_{\delta}\). Also, the middle-third Cantor set is simultaneously \(F_{\sigma}\) and \(G_{\delta}\) \cite{Kechris1995}.\par

\subsection{General Representation}\label{2.2.GR}

At the broadest level, a Cantor set is identified with whatever closed remainder is left in $[0,1]$ after excising a countable family of pairwise disjoint open intervals; no further constraint is placed on how those intervals are chosen, distributed, or scaled. This representation therefore captures the largest class of objects that the term ``Cantor set'' can
reasonably denote on the real line. A precise set-theoretic formulation of this viewpoint was given by Astels \cite{Astels1999}. To begin with,  let $I_{0}=[0,1] $ be the closed unit interval and let  $\{I_{n}^{*}\}_{n\ge1} $ be a (finite or countable) family of pairwise–disjoint open
sub‑intervals contained in $I_{0}$. The general definition of Cantor set $\Gamma$ is given by:
\begin{eqnarray}\label{eq1}
    \Gamma=I_0-\bigcup_{n=1}^{+\infty} I_{n}^{*}.
\end{eqnarray}
Considering closed sets $I_{n}=I_0-I_{n}^{*} (n\in\mathbb{N})$ we have the following equivalent for the most general representation of Cantor set:
\begin{eqnarray}\label{eq2}
    \Gamma=\bigcap_{n=0}^{+\infty} I_{n}.
\end{eqnarray}
This definition covers many non-fractal subsets of the closed unit interval. Furthermore, by definition $\Gamma$ is closed subset of $[0,1]$ and hence a $F_{\sigma}$ set. Also, given that for the Euclidean distance $d_E$ we have $\Gamma=\cap_{n=1}^{+\infty} \Gamma^{(n)}$ where $\Gamma^{(n)}:=\{x\in [0,1] | d_E(x,\Gamma) < \frac{1}{n}\}  (n\in\mathbb{N})$ are open sets, it follows that $\Gamma$ is a $G_{\delta}$ set.\par 

\subsection{Nested Representation}\label{2.3.NR}
The nested representation refines the general one by demanding that the construction proceed in stages: at every step, the surviving set is a finite union of closed intervals, each of which is subdivided in the next step, with diameters tending uniformly to zero. This staged structure---absent from the general representation---is what enables stage-by-stage geometric and measure-theoretic analysis. A formal version of this more structured viewpoint was introduced by Falconer in 2014 \cite{Falconer2014} where a nested structure is imposed on their construction process as follows. Let \(\bigl(s_n\bigr)_{n \ge 0}\) be a strictly increasing sequence of natural numbers with \(s_0 = 1\) and \(s_n \to \infty\).
Construct closed sets \(I_n \subset [0,1]\) for given $n\geq 0$ recursively as follows:

\begin{enumerate}
\item Base step. Set \(I_0 = [0,1]\) and denote its single component by \(I_{0,1}\).
\item Inductive step. Assume \(I_{\,n-1} = \bigcup_{k=1}^{s_{\,n-1}} I_{\,n-1,k}\) is the disjoint union of \(s_{\,n-1}\) closed intervals. Subdivide each \(I_{\,n-1,k}\) into $m_{n,k} \ge 1$ pairwise–disjoint closed sub‑intervals. From the entire collection of these sub‑intervals choose exactly \(s_n\) of them (at least one descendant from every parent), relabel them \(I_{\,n,1},\dots,I_{\,n,s_n}\), and set:
\begin{eqnarray}
I_n &=& \bigcup_{k=1}^{s_n} I_{\,n,k} \subseteq I_{\,n-1}\ \ (n\in\mathbb{N}). 
\end{eqnarray}
 
\end{enumerate}
Assume moreover that the intervals shrink, i.e.
$\max_{1 \le k \le s_n} \lvert I_{\,n,k} \rvert \;\xrightarrow[]{}\; 0.$ Then $I_0 \supset I_1 \supset I_2 \supset \cdots,\text{with}\   I_n = \bigcup_{k=1}^{s_n} I_{\,n,k} \  \text{and}\ \ s_n \uparrow +\infty.$
The associated Cantor‑set limit set is:
\begin{eqnarray}\label{eq3}
\Gamma  &=& \bigcap_{n=0}^{+\infty} I_n.    
\end{eqnarray}
This definition encompasses a broad class of real-line fractals, characterized by their symmetry, Lebesgue measure, and key topological attributes—most notably perfectness and nowhere denseness.\par 

\subsection{Iterated Function System (IFS) Representation}\label{2.4.IFSR}
\noindent
The IFS representation further specializes the nested construction by requiring that the surviving sub-intervals at every stage arise as images of $[0,1]$ under a fixed finite family of affine contractions. The resulting Cantor set is then the unique non-empty compact attractor of the system, and its self-similarity becomes a structural identity rather than an asymptotic property. This more rigid formulation is naturally expressed in the language of iterated-function systems, and appears as a special instance of the nested representation in Falconer \cite{Falconer2014}.  In this case, at each stage of the construction the surviving sub-intervals arise as images of the unit interval under a fixed finite family of affine maps. To see construction,  we recall the relevant terminology as follows. First, we consider the following special maps:

\begin{enumerate}
  \item A map \(T:[0,1]\to[0,1]\) is \emph{affine} if
        \(T(tx_1+(1-t)x_2)=tT(x_1)+(1-t)T(x_2)\) for all
        \(t,x_1,x_2\in[0,1]\).
  \item It is a \emph{contraction} if
        \(\lvert T(x_1)-T(x_2)\rvert\le r\lvert x_1-x_2\rvert\)
        for some \(r\in(0,1)\) and all \(x_1,x_2\in[0,1]\).
\end{enumerate}

Second, every affine contraction can be written \(T(x)=rx+s\) with
\(r,s\in[0,1]\) and \(r<1\) \cite{Edgar2007}; we call
\(r=\operatorname{cont.coeff}(T)\) its \emph{contraction coefficient}.
Third, let \(\mathcal T=\{T_k\}_{k=1}^{K}\;(K\ge2)\) be a finite family of
affine contractions with
\(\max_{1\le k\le K}\operatorname{cont.coeff}(T_k)<1\).
Such a family constitutes an IFS.
For a non-empty compact set \(I\subset[0,1]\) define:
\begin{eqnarray}
T(I) &=&\bigcup_{k=1}^{K} T_k(I).    
\end{eqnarray}
The operator \(T\) is itself a contraction with ratio
\(r=\max_{1\le k\le K}\operatorname{cont.coeff}(T_k)\) \cite{Falconer2014}.

 Hence the sequence $\{T^{\circ n}([0,1])\}_{n=0}^{+\infty}$, where
$T^{\circ n}$ denotes the $n$-fold composition and $T^{\circ 0} = \mathrm{id}$,
converges in the Hausdorff metric to a unique non-empty compact set, called
the \emph{attractor} of the system, and we set the Cantor set as:  \par 
 
\begin{eqnarray}
\Gamma  &=& \bigcap_{n=0}^{+\infty} T^{\circ n}([0,1]).
\end{eqnarray}
Here, each \(T^{\circ n}([0,1])\) is the disjoint union of \(K^n\) closed
intervals of length at most \(r^{\,n}\), so their diameters tend to
$0$ as $n\longrightarrow +\infty.$  Furthermore, the Cantor set  satisfies the
invariance relation:
\begin{eqnarray}
\Gamma  &=& T(\Gamma)\;=\;\bigcup_{k=1}^{K} T_k(\Gamma).    
\end{eqnarray}

Thus, the affine–IFS description furnishes the more specific and informative representation in the nested hierarchy of set-theoretic representations of Cantor sets.\par 

\subsection{q-Ary Expansion}\label{2.5.qAE}
  The $q$--ary expansion is the most rigid level of the hierarchy: every affine contraction in the IFS is required to share the common slope $1/q$, so each retained interval at stage $n$ has length exactly $q^{-n}$ and each point of the limit set admits a digit-string address in $\mathcal{A}^{\mathbb{N}}$. This arithmetization is what makes explicit closed-form descriptions possible. The $q-$Ary expansion of Cantor sets is indeed an special case of their IFS representation when the involved affine maps have equal slopes. In details, let $q\in\mathbb{N}+2,$ and $ A\subset \{0,1,\cdots,q-1\}, |A|=K (2\leq K \leq q ).$ For every digit $a\in A$ define the $a^{th}$ strict contraction $T_a$ with slope $\frac{1}{q}$ and y-intercept $\frac{a}{q}:$
\begin{eqnarray}
T_a&:& [0,1]\longrightarrow [0,1] \nonumber\\
T_a(x)&=&\frac{1}{q} x + \frac{a}{q}.
\end{eqnarray}
Now, the IFS representation of the Cantor set $\Gamma$ is given by:
\begin{eqnarray}
\Gamma&=& \bigcap_{n=0}^{+\infty} I_n: 
I_0=[0,1], I_n= \bigcup_{a\in A}  T_a(I_{n-1}),\ \   (n\in\mathbb{N}).    
\end{eqnarray}
This special case has the following properties: 
\begin{enumerate}
\item Each component of $I_n$ is a closed interval of length $q^{-n}.$
\item There are $K^n$ such components with one for each finite digit block $\vec{a}=(a_1,\cdots,a_n): a_k\in A (1\leq k\leq n).$
\item The system is nested:  $I_n \supset I_{n-1} (n\in\mathbb{N}).$
\end{enumerate}
A straightforward verification by mathematical induction shows that:
\begin{eqnarray}
I_n&=& T_{a_1} \circ \cdots \circ T_{a_n} ([0,1])= [\sum_{k=1}^{n}\frac{a_k}{q^k}, \sum_{k=1}^{n}\frac{a_k}{q^k}+\frac{1}{q^n}], \ \ (n\in\mathbb{N}). 
\end{eqnarray}
Also, the map  $F_n=T_{a_1} \circ \cdots \circ T_{a_n}$  is $q^{-n}$-Lipschitz continuous $(n\in\mathbb{N})$ \cite{Searcoid2006}. Consequently, $x\in \Gamma $ if and only if there is a unique (up to the standard endpoint ambiguity) infinite digit string $\vec{a}^{*}=(a_1,a_2,a_3,\cdots)\in A^{\mathbb{N}}$ such that $x=\sum_{n=1}^{+\infty} \frac{a_n}{q^n}.$ Hence, we have the $q-$ ary expansion of the Cantor set:
\begin{eqnarray}
\Gamma&=& \{x\in[0,1]| x= \sum_{n=1}^{+\infty} \frac{a_k}{q^n}, a_n\in A \}.    
\end{eqnarray}
\textbf{Notation:} Henceforth, we denote $\Gamma(\alpha, K)$ middle $\alpha-$ Cantor set defined by $K$ affine maps.  \par 

In the upcoming subsections, we study, explicit and recursive representation of three families of Cantor sets given their associated $q-$ary expansions. \par 
  
\subsection{Hierarchical Relationships and Practical Differences}
\label{2.6.HRPD}

The four representations introduced in Sections~\ref{2.2.GR}-\ref{2.5.qAE} form a strict descending chain of structural specificity:
\begin{equation}
\text{General} \;\supsetneq\; \text{Nested} \;\supsetneq\; \text{IFS}
\;\supsetneq\; q\text{--ary}.
\label{eq:hierarchy}
\end{equation}
Each inclusion is proper, and at each step a concrete structural commitment is added to the preceding level:
\begin{itemize}
\item \emph{General $\to$ Nested.} The arbitrary countable family of
    deleted open intervals is reorganized into a stagewise construction in which the surviving set at each stage is a finite union of closed intervals with diameters tending to zero. The strictness of the inclusion follows from the existence of closed nowhere-dense subsets of $[0,1]$ that admit a representation of the form~\eqref{eq1} but not a stagewise nested decomposition with $s_n \uparrow \infty$
    (see~\cite{Astels1999}).
    \item \emph{Nested $\to$ IFS.} The stagewise construction is required to be generated by a fixed finite family of affine contractions, so that the limit set satisfies the self-similarity identity $\Gamma = \bigcup_{k=1}^K T_k(\Gamma)$. Not every nested Cantor set is an IFS attractor: nested constructions in which the contraction ratios vary across stages or across components fail this identity
    (cf.~\cite{Hutchinson1981, Falconer2014}).
    \item \emph{IFS $\to$ $q$--ary.} All affine contractions are forced to share the common slope $1/q$ for some integer $q \geq 2$. IFS attractors with unequal contraction ratios---e.g.\ the standard examples with $r_1 \neq r_2$ used to realize prescribed Hausdorff dimensions---are excluded.
\end{itemize}
The price of descending the hierarchy is generality; the reward is
descriptive power. Table~\ref{tab:hierarchy} summarizes the practical
trade-offs.

\begin{table}[H]
\caption{Practical comparison of the four set-theoretic representations of Cantor Sets on the real line.}
\centering
\small
\begin{tabular}{@{}p{0.4cm}p{2.8cm}p{2.8cm}p{2.8cm}p{2.8cm}p{2.8cm}@{}}
\toprule
\#&\textbf{Item} & \textbf{General} & \textbf{Nested} & \textbf{IFS} & \textbf{$q$--ary} \\
\midrule
1&Operative Parameter
& Family $\{I_n^*\}$ of deleted open intervals
& Sequence $(s_n)$ and per-stage subdivisions
& Affine maps $\{T_k\}_{k=1}^{K}$
& Base $q$ and digit set $\mathcal{A} \subseteq \{0,\dots,q-1\}$ \\[2pt]
 & & & & & \\[2pt]
2&Stage-$n$ Structure
& Not stage-defined in general
& $s_n$ closed intervals of varying length
& $K^n$ closed intervals, length $\leq r^n$
& $K^n$ closed intervals of length exactly $q^{-n}$ \\[2pt]
 & & & & & \\[2pt]
3&Self-Similarity
& None required
& None required
& Exact: $\Gamma = \bigcup_k T_k(\Gamma)$
& Exact, with uniform ratio $1/q$ \\[2pt]
 & & & & & \\[2pt]
4&Closed-form Points
& Not available in general
& Not available in general
& Available via fixed-point address
& Available via digit expansion $\sum a_n q^{-n}$ \\[2pt]
 & & & & & \\[2pt]
5&Typical Use
& Existence and Borel-class arguments
& Stagewise measure / dimension analysis
& Self-similarity, dimension theorems
& Explicit formulas, arithmetic structure \\
\bottomrule
\end{tabular}
\label{tab:hierarchy}
\end{table}

The three families treated in Section~\ref{2.STF} illustrate the lower end of this hierarchy concretely: $\Gamma_1$ and $\Gamma_2$ are presented through their $q$--ary expansions and admit closed-form gap descriptions, while $\Gamma_3(p/q,2)$, although IFS-generated, requires the recursive endpoint treatment of Theorem~\ref{thm3} because its retained-interval lengths do not collapse to a single power of a fixed base.\par

\section{Descriptive Formulas for Cantor Sets}\label{3.DFCS}

\subsection{Thin Family \#1}\label{3.1.TF1}
The earliest set-theoretic explicit description of Cantor sets appears in 2006 \cite{Soltanifar2006a}, where the author treats the symmetric thin family generated by at least two affine maps in an iterated-function system. This construction yields a fully constructive proof of the Hausdorff-Dimension Theorem  (which asserts that for every $s>0$ there exists a continuum of compact subsets of $\mathbb{R}^n\quad (n\geq \lceil s \rceil)$ with Hausdorff dimension of exactly $s$), thereby establishing the existence of the corresponding fractals for any prescribed Hausdorff dimension. Subsequent work extends the argument to broader deterministic and stochastic settings \cite{Soltanifar2021,Soltanifar2022}.\par   

\begin{thm}\label{thm1}
Let  $\Gamma_1(\frac{1}{q}, \frac{q+1}{2}):=\{ x\in [0,1] | x=\sum_{n=1}^{+\infty}\frac{a_n}{q^n} \ , \ a_n=0,2,4,\cdots,q-1 \}$ where  $q\in 2\mathbb{N}+1.$ Then:
\begin{eqnarray}
\Gamma_1(\frac{1}{q},\frac{q+1}{2}) &=& [0,1] - \bigcup_{n=1}^{+\infty} \bigcup_{k=0}^{q^{n-1}-1} \bigcup_{r=1}^{\frac{q-1}{2}} (\frac{qk+(2r-1)}{q^n}, \frac{qk+2r}{q^n}): \ \   (q\in 2\mathbb{N}+1).  
\end{eqnarray}
\end{thm}
In particular, for the case $q=3$ we obtain the explicit formulae for the middle-third Cantor set $C=\Gamma_1(\frac{1}{3},2)$ as in the second formulae in Corollary \ref{cormain}.\par 

\subsection{Thin Family \#2}\label{3.2.TF2} 
In the same year,  the explicit set-theoretic description of a second family of Cantor sets was presented  \cite{Soltanifar2006b}. As with the earlier construction, the sets obtained are thin and symmetric; however, their underlying iterated-function system involves only two affine contractions, in contrast to  generally larger number used for the first family.\par  

\begin{thm}\label{thm2}
Let  $\Gamma_2(\frac{q-2}{q}, 2):=\{ x\in [0,1] | x=\sum_{n=1}^{+\infty}\frac{a_n}{q^n} \ , \ a_n=0,q-1 \}$ where  $q\in \mathbb{N}+2.$ Then:
\begin{eqnarray}
\Gamma_2(\frac{q-2}{q},2) &=& [0,1] - \bigcup_{n=1}^{+\infty} \bigcup_{k=0}^{q^{n-1}-1} (\frac{qk+1}{q^n}, \frac{qk+2}{q^n}): \ \   (q\in  \mathbb{N}+2).  
\end{eqnarray}
\end{thm}
Similar to the Thin Family $\Gamma_1$, for the case $q=3$ we obtain the explicit formulae for the middle-third Cantor set $C=\Gamma_2(\frac{1}{3},2)$ as in the second formulae in Corollary \ref{cormain}.\par 

\subsection{Augmented Thick Family \#1}\label{3.2.ATF}  
In both preceding cases we provided closed-form descriptions for two families of measure-zero Cantor sets. This naturally raises the question: Does an analogous explicit representation exist for Cantor sets of positive Lebesgue measure? Current evidence suggests a negative answer; no such closed formula is known at the time of writing. Nevertheless, a constant-recursive representation remains feasible \cite{Elaydi2005}, and the requisite background together with interim results is presented in Appendix A for interested readers.\par 

\begin{figure}[H] 
\centering 
\includegraphics[clip,width=0.75\columnwidth]{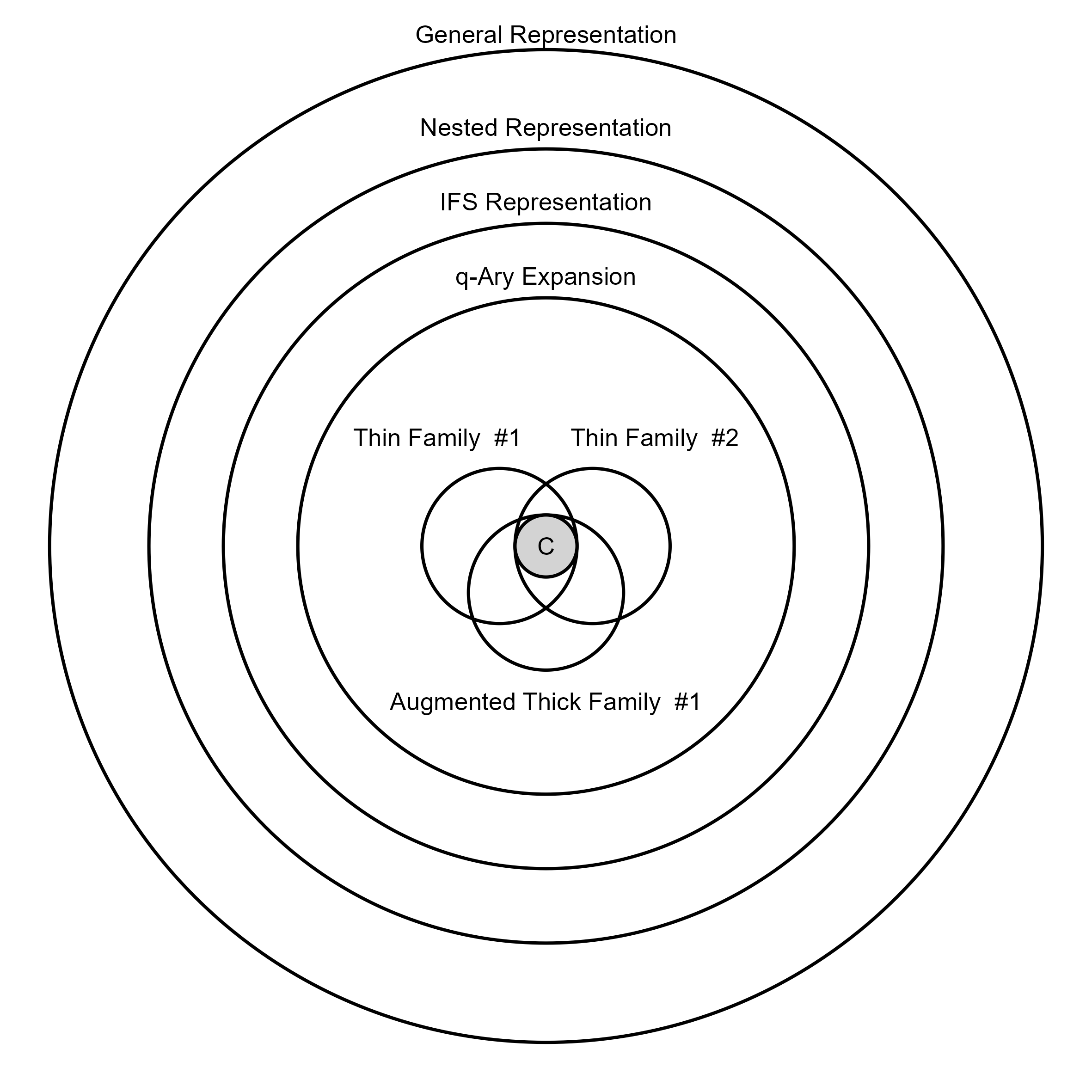}
\caption{Diagram of descriptive formulas for the Cantor sets from set-theoretic perspective: C refers to the middle-third Cantor set.\label{fig1}}
\end{figure}

\begin{thm}\label{thm3}
Let  $\Gamma_3(\frac{p}{q}, 2):=\{ x\in [0,1] | x=\sum_{n=1}^{+\infty}\frac{a_n}{(2q)^n} \ , \ a_n=0,1,\cdots,q-p-1,q+p,\cdots,2q-1,\}$ where  $p,q\in \mathbb{N}, (p,q)=1, 0<\frac{p}{q}\leq\frac{1}{3}.$ Then:
\begin{eqnarray}
\Gamma_3(\frac{p}{q},2) &=& 
\bigcap_{n=0}^{+\infty}\bigcup_{k=1}^{2^n}  [ a_{n,k}(\frac{p}{q}), b_{n,k}(\frac{p}{q})]:\\
&& a_{0,1}(\frac{p}{q})=0;\nonumber\\
&& a_{n,k}(\frac{p}{q})=1_{odd}(k)\Big(a_{n-1,\frac{k+1}{2}}(\frac{p}{q})\Big) + 1_{even}(k)\Big(a_{n-1,\frac{k}{2}}(\frac{p}{q})+\frac{(1-3\frac{p}{q})+(1-\frac{p}{q})(2\frac{p}{q})^n}{(1-2\frac{p}{q})2^n}\Big);\nonumber \\
&& b_{0,1}(\frac{p}{q})=1;\nonumber\\
&& b_{n,k}(\frac{p}{q})=1_{odd}(k)\Big(b_{n-1,\frac{k+1}{2}}(\frac{p}{q})-\frac{(1-3\frac{p}{q})+(1-\frac{p}{q})(2\frac{p}{q})^n}{(1-2\frac{p}{q})2^n}\Big) + 1_{even}(k)\Big(b_{n-1,\frac{k}{2}}(\frac{p}{q})\Big): \nonumber\\
&& 1\leq k\leq 2^n, n\in\mathbb{N}\nonumber.
\end{eqnarray}
\end{thm}
The  Computer software code to calculate the  end points $a_{n,k}(.), b_{n,k}(.)$ are given in  Appendix B with associated examples for $\alpha_1=\frac{1}{3},$ and $\alpha_2=\frac{1}{4}$ \cite{RCoreTeam2022}. It is noteworthy to mention that the middle-third Cantor set belongs to all three discussed families. Figure \ref{fig1} summarizes our above discussions in terms of hierarchy of descriptive representations.\par 

Finally, for the case $p=1,q=3$ we obtain the first order recursive formulae for the middle-third Cantor set $C=\Gamma_3(\frac{1}{3},2)$  as in the third formulae in the following Corollary \ref{cormain}: 
\begin{cor}\label{cormain}
The middle-third Cantor set  has three equivalent set theoretic  descriptive representations:
\begin{eqnarray}
C&=& \{ x\in [0,1] | x=\sum_{n=1}^{+\infty}\frac{a_n}{3^n} \ , \ a_n=0,2 \}\\
 &=& [0,1] - \bigcup_{n=1}^{+\infty} \bigcup_{k=0}^{3^{n-1}-1} (\frac{3k+1}{3^n}, \frac{3k+2}{3^n})\\
 &=& \bigcap_{n=0}^{+\infty}\bigcup_{k=1}^{2^n}  [a_{n,k}, b_{n,k}]: \\
 && a_{0,1}=0; a_{n,k}=1_{odd}(k) \Big(a_{n-1,\frac{k+1}{2}}\Big)+1_{even}(k)\Big(a_{n-1,\frac{k}{2}}+\frac{2}{3^n}\Big); \nonumber\\
 &&
b_{0,1}=1; b_{n,k}=1_{odd}(k) \Big(b_{n-1,\frac{k+1}{2}}-\frac{2}{3^n}\Big)+1_{even}(k)\Big(b_{n-1,\frac{k}{2}}\Big); \nonumber\\
&& 1\leq k\leq 2^n, n\in\mathbb{N}\nonumber. 
\end{eqnarray}
\end{cor}

\section{Discussion}\label{4.D}

\subsection{Summary}\label{4.1SC}
  The four descriptions surveyed here---general, nested, IFS, and $q$--ary---form a strict hierarchy in which each level imposes additional structure on its predecessor, and the three families $\Gamma_1$, $\Gamma_2$, $\Gamma_3$ populate this hierarchy at successively finer levels of specificity. The recursive endpoint formula for $\Gamma_3(p/q,2)$ closes a gap in the positive-measure case: where the thin families admit closed-form gap descriptions, the augmented thick family had previously resisted such an explicit treatment, and the recursion presented in Theorem~\ref{thm3}
now makes every retained interval at every finite stage exactly computable. The middle-third Cantor set, sitting in the intersection of all three families, illustrates how a single canonical object can serve as a stress test for the entire descriptive framework. The remainder of this section explores three consequences of this synthesis: a continuum-theoretic reading of $\Gamma_3$ via embeddings into standard continua, a measure-theoretic construction of canonical self-similar measures supported on $\Gamma_3(\alpha)$,
and an outlook tying the framework to the open problems of Section~\ref{4.2FWOP}.   \par 

\paragraph{Continuum-theoretic perspective.}
Although Cantor sets are totally disconnected, they occur naturally inside many \emph{continua} (compact connected metric spaces) as closed, perfect, nowhere-dense subsets~\cite{Nadler1992}. The explicit endpoint recursion for $\Gamma_3$ supplies fine control of gap geometry at every stage. This facilitates: (i) explicit embeddings of $\Gamma_3(\alpha)$ into standard continua (e.g.,  as closed subsets of \emph{arcs} (homeomorphic images of $[0,1]$), \emph{dendrites} (locally connected continua containing no simple closed curve), or \emph{Peano continua} (locally connected continua, equivalently continuous images of $[0,1]$) ) with prescribed gap patterns; (ii) constructive “fill-in’’ procedures that yield continua whose complements are modeled by a chosen $\Gamma_3(\alpha)$; and (iii) numerical investigation of intersection phenomena (useful in studying when two embedded Cantor-type sets intersect), where computable gap data and thickness proxies are central. In this sense the paper’s formulas act as a bridge, turning qualitative continuum-theory questions into explicit, checkable computations.\par 

\paragraph{Measure-theoretic perspective (a canonical $\Gamma_3$ measure).}
For each $\Gamma_3(\alpha)$ we can \emph{canonically} equip the set with a Borel probability measure that generalizes the Cantor–Lebesgue measure \cite{Rosenthal2006}. At stage $n$ there are $q^n$ retained intervals; assign either uniform weights $q^{-n}$ or, more generally, a fixed weight vector $w=(w_1,\dots,w_q)$ with $\sum_{i=1}^{n} w_i=1$, propagated multiplicatively along the recursion. This produces a self-similar probability measure $\mu_{\alpha,w}$ supported on $\Gamma_3(\alpha)$ whose cumulative distribution function is a Devil’s–staircase generalization of the Cantor function: it is continuous, nondecreasing, constant on deleted gaps, and increases exactly on $\Gamma_3(\alpha)$. The explicit endpoints $(a_{n,k},b_{n,k})$ allow fast, exact evaluation of cylinder-set masses and numerical plots of the associated distribution function. This construction is immediate from the framework developed here and is implemented by minor extensions of Appendix~B.

\paragraph{Outlook.}
Because the geometry of $\Gamma_3$ is explicit at every finite stage, the framework supports systematic computation of gap-length statistics, empirical thickness proxies, and intersection behavior under translation—tools that are relevant both to continuum-theoretic embeddings and to questions in additive/fractal geometry. We conclude with open problems on asymmetry, irrational endpoints, variable-gap deletions, the role of symmetry, and common intersections; these problems are now posed in a setting where algorithmic experiments and exact finite-stage calculations are straightforward.

\subsection{Future Work and Open Problems}\label{4.2FWOP}

The preceding survey has shown that three natural families of deterministic
Cantor sets---two of measure zero with explicit descriptions and one of
non-negative measure defined recursively---share a common symmetric
architecture and possess rational construction end-points; indeed, the
classical middle-third Cantor set lies in the intersection of all three.
These parallel features highlight several directions in which the present
set-theoretic framework remains incomplete and invite further research.
We organize the open problems into three thematic clusters reflecting
their logical interdependence: \emph{structural generalization}
(Theme~A), \emph{variable-scale constructions} (Theme~B), and
\emph{interaction phenomena} (Theme~C), with the last depending on
progress in the first two.

\subsubsection*{Theme A. Structural Generalization}

The descriptive formulas of Section~\ref{3.DFCS} rely on two structural
restrictions---bilateral symmetry of the deletion pattern and rationality
of the retained-interval endpoints. The following two problems ask whether
either restriction can be dropped without losing explicitness.

\begin{itemize}
    \item[OP1.] \emph{Asymmetric constructions.} Devise explicit or
    recursive set-theoretic formulas for Cantor sets generated by
    iterated-function systems with $K \geq 2$ maps that do not preserve
    bilateral symmetry.
    \item[OP2.] \emph{Irrational endpoints.} Find descriptive formulas
    when the retained intervals have irrational endpoints, again
    allowing $K \geq 2$ IFS components.
\end{itemize}

\paragraph{Rationale and impact.}
Symmetric, rationally-located Cantor sets form a measure-zero subfamily
within the space of all IFS attractors, yet they are essentially the only
ones for which closed or recursive descriptive formulas are currently
known. A resolution of OP1 would extend the techniques of Section~\ref{3.DFCS}
to the asymmetric self-similar sets that arise as hyperbolic invariant
sets in low-dimensional dynamics~\cite{Palis1993, Moreira2001}, where
asymmetry is the rule rather than the exception. A resolution of OP2
would connect the descriptive framework to the algebraic theory of
self-similar measures with non-rational contraction ratios, where deep
recent work by Hochman~\cite{Hochman2014} and Shmerkin~\cite{Shmerkin2019}
has clarified the dimension theory but explicit set-theoretic formulas
remain absent.

\paragraph{Tools and feasibility.}
The endpoint recursion of Theorem~\ref{thm3} is candidate machinery for
OP1: removing the symmetry assumption requires tracking left- and
right-endpoints under independent contraction ratios $r_L \neq r_R$,
which extends the recursion in Lemmas~4.1--4.3 (Appendix~A) at the cost
of a more elaborate index scheme. OP2 is genuinely harder, since the
arithmetic identity $p/q = 2p/(2q)$ that drives the $2q$-ary expansion
in Appendix~A has no analogue for irrational ratios; symbolic-dynamic
coordinates~\cite{Lind1995} or continued-fraction expansions of the
endpoints offer alternative descriptive routes that have not yet been
systematically explored.

\subsubsection*{Theme B. Variable-Scale Constructions}

The middle-$\alpha$ family $\Gamma_3(\alpha,2)$ removes a gap of length
$\alpha^n$ at stage $n$, so the deletion sequence is geometric. The
following two problems concern non-geometric deletion schemes and the
role of symmetry within them.

\begin{itemize}
    \item[OP3.] \emph{Variable-gap deletions.} Given a sequence
    $(\alpha_n)_{n \geq 1}$ with $0 < \alpha_n < 1$ satisfying
    $1 - \sum_{n=1}^{\infty} 2^{n-1} \alpha_n \geq 0$, provide a
    recursive descriptive formula for the Cantor set obtained by removing
    at stage $n$ the open middle interval of length $\alpha_n$.
    \item[OP4.] \emph{Role of symmetry.} Determine how the solution of
    OP3 varies when bilateral symmetry of the deletion is enforced or
    relaxed.
\end{itemize}

\paragraph{Rationale and impact.}
Variable-scale constructions interpolate between the self-similar Cantor
sets of Section ~\ref{2.STF} and the broader nested representations of
(\ref{eq3}). Smith--Volterra--Cantor sets and their generalizations belong
to this class~\cite{Valline2013}, and their Lebesgue measure can be tuned
arbitrarily within $[0,1)$ by choice of $(\alpha_n)$. A resolution of
OP3 would yield exact endpoint formulas for these positive-measure
constructions, parallel to what Theorem ~\ref{thm3} achieves in the
self-similar case, and would supply concrete examples for testing
fractal-uncertainty principles in the variable-scale setting~\cite{Laba2009}.
OP4 is significant because asymmetric variable-gap constructions exhibit
qualitatively different harmonic-analytic behavior---in particular, the
Fourier decay of the associated self-similar measures depends sensitively
on whether deletions are symmetric.

\paragraph{Tools and feasibility.}
The recursion of Lemma~4.3 generalizes directly to OP3 by replacing
$\alpha^n$ with $\alpha_n$ in the increment $\Delta_{n,k}$; the closed
form of $\Delta_{n,k}$ no longer collapses but remains explicit as a
finite sum at each stage. OP4 admits a tractable first case (left/right
asymmetry only, with all other parameters symmetric), which would be a
natural starting point.

\subsubsection*{Theme C. Interaction Phenomena}

Themes A and B describe families of Cantor sets indexed by structural
parameters (asymmetry, irrationality, deletion sequence). The next
problem asks how families from these themes interact.

\begin{itemize}
    \item[OP5.] \emph{Common ground.} Establish conditions under which
    the families arising in OP1--OP4 admit a non-empty mutual
    intersection, and characterize the intersection when it exists.
\end{itemize}

\paragraph{Rationale and impact.}
The middle-third Cantor set is, in the present survey, the unique
explicit point of intersection across the three families
$\Gamma_1, \Gamma_2, \Gamma_3$. Whether analogous canonical objects
exist across the broader families of OP1--OP4 is a question of both
descriptive and dynamical interest. Common intersections are the
deterministic shadow of the Newhouse phenomenon and the Palis
conjecture~\cite{Newhouse1979, Palis1993, Moreira2001}: when do two
Cantor sets meet, and when does the meeting persist under perturbation?
A resolution of OP5 would supply concrete test cases for these
conjectures and would bridge the explicit set-theoretic framework of
this survey to the geometric measure-theoretic literature on
intersections~\cite{Mattila1995}.

\paragraph{Tools and feasibility.}
The thickness theory of Astels~\cite{Astels1999, Astels2000} provides
sufficient conditions for non-empty intersection of two Cantor sets in
terms of their gap structures, and the explicit endpoint data supplied
by Theorem~\ref{thm3} (and its Theme-A and Theme-B extensions) make
these thickness computations exact rather than asymptotic.
Moreira--Yoccoz~\cite{Moreira2001} provides the deeper machinery for
\emph{stable} intersections under perturbation, which becomes relevant
once OP1--OP4 have been resolved.

\paragraph{}
Resolving these questions would extend the taxonomy of Cantor sets
beyond the symmetric, rational, geometrically-scaled cases treated here
and clarify how set-theoretic representations interact with geometric
asymmetry, irrational numbers, variable deletion schemes, and the
intersection phenomena that link them. The thematic grouping above is
intended to make explicit the logical dependence---Theme~C builds on
Themes~A and B---and to facilitate parallel progress on the structural
and variable-scale fronts.

\section*{Conclusion}\label{5.C}
This survey brings four set-theoretic perspectives—general, nested, IFS, and 
$q-$ary—under one framework and demonstrates how they meet in three representative Cantor-set families containing the middle-third classic. The synthesis clarifies which descriptive formulas are now complete and which remain missing, especially for asymmetric, irrational-endpoint, and variable-gap constructions. By distilling these gaps into a concise list of open problems, the paper sets a clear agenda for future work in descriptive set theory and fractal analysis on the real line.\par

\subsection*{Funding} This research received no external funding.
\subsection*{Institutional Review Board Statement} Not Applicable.
\subsection*{Informed Consent Statement} Not Applicable.
\subsection*{Data Availability Statement} Not Applicable.
\subsection*{Acknowledgments} The author is grateful to the reviewers for their constructive comments, which significantly improved the presentation of this manuscript.
\subsection*{Conflicts of Interest} The author  declares no conflicts of interest.
\subsection*{Abbreviations}
The following abbreviations are used in this manuscript:\newline
IFS: Iteration Function System; MC: Main Contribution; OP: Open Problem; SVC: Smith Volterra Cantor
\section*{Appendix}
\subsection*{Appendix.A. Interim Results}\label{Appendix.A}
Let $\Gamma_3(\alpha,2)$ where $0<\alpha\leq\frac{1}{3}$ be the middle $-\alpha$ Cantor set defined by two affine maps. Here, its nested representation is given by $\Gamma_3(\alpha,2)=\cap_{n=0}^{+\infty}I_n(\alpha)$ where $I_0(\alpha)=[0,1]$ and $I_n(\alpha)$ is defined recursively  by  removing $2^{n-1}$ middle open intervals each of the length $\alpha^n$ from  $I_{n-1}(\alpha)$ $(n\in\mathbb{N}).$ We have the following interim results:
\begin{enumerate}
\item \textbf{Family Name} By construction, the Lebesgue measure is given by $\lambda(\Gamma_3(\alpha,2))=1-\sum_{n=1}^{+\infty}2^{n-1}\alpha^n=\frac{1-3\alpha}{1-2\alpha} (0<\alpha\leq\frac{1}{3}).$ Hence, we augment the "Thick Family" with $0<\alpha<\frac{1}{3}$ with the "Thin case" $\alpha=\frac{1}{3};$ and, call it the "Augmented Thick Family" with $0<\alpha\leq\frac{1}{3}.$ 

\item \textbf{2q-Ary Expansion} When $0<\alpha=\frac{p}{q}\leq\frac{1}{3}, p,q\in\mathbb{N}, (p,q)=1$, by above construction, the equality $\frac{p}{q}=\frac{2p}{2q}$, and mathematical induction, a straightforward verification shows the following 2q-Ary expansion: $$\Gamma_3(\alpha,2)=\{ x\in [0,1] | x=\sum_{n=1}^{+\infty}\frac{a_n}{(2q)^n} \ , \ a_n=0,1,\cdots,q-p-1,q+p,\cdots,2q-1\}.$$
This representation facilitates the n-Aray expansion for cases with even $q.$ Furthermore, the middle-third Cantor set has the  secondary 6-Ary expansion as $C=\{ x\in [0,1] | x=\sum_{n=1}^{+\infty}\frac{a_n}{6^n} \ , \ a_n=0,1,4,5\}.$   

\item \textbf{Proof of Theorem \ref{thm3}} We prove the Theorem in  steps (a)-(e). 
\begin{enumerate}
\item 
\textbf{Notation}
Given above recursive construction it follows that,   $\Gamma_3(\alpha,2)=\cap_{n=0}^{+\infty}I_n(\alpha)$ where $I_n(\alpha)=\cup_{k=1}^{2^n} I_{n,k}(\alpha)$ and $I_{n,k}(\alpha)=[a_{n,k}(\alpha), b_{n,k}(\alpha)], (1\leq k\leq 2^n, n\in\mathbb{N}_{0}).$
\item
\begin{lem}\label{lemA.1}
$\lambda(I_{n,k}(\alpha))=\frac{(2\alpha)^{n+1}+2(1-3\alpha)}{2^{n+1}(1-2\alpha)}: \ \ (1\leq k\leq 2^n, n\in\mathbb{N}_{0}).$
\end{lem}
\textbf{Proof.}  The case for $n=0$ is trivial. Let $n\in\mathbb{N}$ and $J_{n,k}(\alpha)$ be the middle open interval of length $\alpha^n$ removed from $I_{n-1,k}(\alpha)$ in the nested construction process i.e. $I_{n-1,k}(\alpha)=I_{n,2k-1}(\alpha)\cup J_{n,k}(\alpha)\cup I_{n,2k}(\alpha).$ Then, given, $\lambda(I_{n,2k-1}(\alpha))=\lambda(I_{n,2k}(\alpha))$ we have: $\lambda(I_{n-1,k}(\alpha))=2\lambda(I_{n,k}(\alpha))+\alpha^n.$ Consequently, we have a first order recursive relation:
$$\lambda(I_{n,k}(\alpha))=\frac{1}{2} \lambda(I_{n-1,k}(\alpha))-\frac{\alpha^n}{2}, \lambda(I_{0,1}(\alpha))=1.$$
Finally, the assertion follows by mathematical induction on $n\in\mathbb{N}$. $\Box$
\item 
\begin{lem}\label{lemA.2}
Let $\Delta_{n,k}(\alpha)=\lambda(I_{n,k}(\alpha))+\alpha^n \ \ (1\leq k\leq 2^n, n\in\mathbb{N}_{0}).$ Then, $$\Delta_{n,k}(\alpha)=\frac{(1-3\alpha)+(1-\alpha)(2\alpha)^n}{(1-2\alpha)2^n} \ \ (1\leq k\leq 2^n, n\in\mathbb{N}_{0}).$$   
\end{lem}
\textbf{Proof.} This is straightforward result from Lemma \ref{lemA.1}. $\Box$
\item 
\begin{lem}\label{lemA.3}
Let $a_{n,k}(\alpha), b_{n,k}(\alpha), (1\leq k\leq 2^n , n\in\mathbb{N}_0)$ be defined as in above recursive construction. Then:
\begin{eqnarray}
a_{0,1}(\alpha)&=&0, \nonumber\\
a_{n,k}(\alpha)&=&1_{odd}(k)\Big(a_{n-1,\frac{k+1}{2}}(\alpha) \Big)+1_{even}(k)\Big( a_{n-1,\frac{k}{2}}(\alpha)+\Delta_{n,k}(\alpha)\Big), \nonumber\\
b_{0,1}(\alpha)&=&1, \nonumber\\
b_{n,k}(\alpha)&=&1_{odd}(k)\Big(b_{n-1,\frac{k+1}{2}}(\alpha)-\Delta_{n,k}(\alpha)\Big)+1_{even}(k)\Big( b_{n-1,\frac{k}{2}}(\alpha) \Big),\nonumber\\
&&(1\leq k\leq 2^n , n\in\mathbb{N}_0). \nonumber 
\end{eqnarray}
\end{lem}
\textbf{Proof.} This is straightforward result from  definition of $\Delta_{n,k}(\alpha)$ in Lemma \ref{lemA.2}, and comparing the end points of closed intervals in the following equations:
\begin{eqnarray}
[a_{n-1,k}(\alpha),b_{n-1,k}(\alpha)]&=&I_{n-1,k}(\alpha)=I_{n,2k-1}(\alpha)\cup J_{n,k}(\alpha)\cup I_{n,2k}(\alpha)\nonumber\\
&=& [a_{n,2k-1}(\alpha),b_{n,2k-1}(\alpha)]
\cup
(b_{n,2k-1}(\alpha),a_{n,2k}(\alpha))
\cup
[a_{n,2k}(\alpha),b_{n,2k}(\alpha)],\nonumber\\
&&(1\leq k\leq 2^n , n\in\mathbb{N}_0).\nonumber
\end{eqnarray}
$\Box$
\item \textbf{Completion} Finally, by considering $0<\alpha=\frac{p}{q}\leq\frac{1}{3}$ in Lemma \ref{lemA.2} and substituting corresponding $\Delta_{n,k}$ in the equations presented in Lemma \ref{lemA.3}, the proof for the Theorem \ref{thm3} is completed. $\Box$   
\end{enumerate}
\end{enumerate}
\subsection*{Appendix.B. R Software Code}\label{Appendix.B}
The following function in R software code produces the end points of the intervals of the recursive representation of "Augmented Thick Family" $\Gamma_3(\frac{p}{q},2)$ for given pair $(n,k): 1\leq k\leq 2^n, n\in\mathbb{N}.$ It presents two examples of calculations one for the thin Cantor set and another for the thick Cantor set:

\begin{verbatim}
## ========================================================================================
##  Recursive evaluator for lower endpoint a_{n,k}(p/q)
## ========================================================================================
a_nk <- local({
  ## simple memorization (saves repeated work)
  cache <- new.env(parent = emptyenv())
  f <- function(n, k, p, q) {
    if (n == 0L) {                     # base level
      if (k != 1L) stop("For n = 0 you must take k = 1")
      return(0)
    }
    if (k < 1 || k > 2^n)
      stop("k must satisfy 1 <= k <= 2^n")
    ## build a unique key for memorization
    key <- paste(n, k, p, q, sep = "-")
    if (exists(key, envir = cache, inherits = FALSE))
      return(cache[[key]])
    r <- p / q                         # the ratio p/q
    ## recursive step
    if (k %% 2L == 1L) {               # k is odd
      val <- f(n - 1L, (k + 1L) / 2L, p, q)
    } else {                           # k is even
      prev <- f(n - 1L, k / 2L, p, q)
      term <- ((1 - 3 * r) + (1 - r) * (2 * r)^n) /((1 - 2 * r) * 2^n)
      val <- prev + term
    }
    cache[[key]] <- val                # store and return
    val
  }
  f
})

## ========================================================================================
##  Recursive evaluator for upper endpoint  b_{n,k}(p/q)
## ========================================================================================

b_nk <- local({
  cache <- new.env(parent = emptyenv())  #memorization
  f <- function(n, k, p, q) {
    if (n == 0L) {                  # base case
      if (k != 1L) stop("For n = 0 you must take k = 1")
      return(1)
    }
    if (k < 1L || k > 2L^n)
      stop("k must satisfy 1 <= k <= 2^n")
    key <- paste(n, k, p, q, sep = "-")
    if (exists(key, envir = cache, inherits = FALSE))
      return(cache[[key]])
    r    <- p / q
    term <- ((1 - 3 * r) + (1 - r) * (2 * r)^n) /((1 - 2 * r) * 2^n)
    ## recursive step
    if (k %% 2L) {                               # k  is odd
      val <- f(n - 1L, (k + 1L) %/% 2L, p, q) - term
    } else {                                     # k is even
      val <- f(n - 1L, k %/% 2L, p, q)
    }
    cache[[key]] <- val
    val
  }
  f
})

## ========================================================================================
##  Intervals Endpoints Calculator  
## ========================================================================================
EndpointCalculator <- function(n, p, q) {
  ## ---------- base‑level endpoints ----------
  a0 <- a_nk(0L, 1L, p, q)
  b0 <- b_nk(0L, 1L, p, q)
  ## keep the original “too‑large” denominator
  denom <- q^(n + 1L)
  ## ---------- level‑n endpoints ----------
  res1 <- sapply(1:2^n, function(k) a_nk(n, k, p, q))
  res2 <- sapply(1:2^n, function(k) b_nk(n, k, p, q))
  num1 <- c(a0 * denom, round(res1 * denom))
  num2 <- c(b0 * denom, round(res2 * denom))
  ## ---------- LABELS (divide out one factor of q) ----------
  display_denom <- denom / q   #  q^(n+1) / q  ->  q^n
  col_labs       <- paste0("k=", 0:2^n)
  fractions1     <- paste0(num1 / q, "/", display_denom)
  fractions2     <- paste0(num2 / q, "/", display_denom)
  ## ---------- NEW ROW ----------
  combo <- paste0("[", fractions1, ",", fractions2, "]")
  out1   <- rbind("a_(n,k)(α)"              = fractions1,
                 "b_(n,k)(α)"              = fractions2,
                 "[a_(n,k)(α),b_(n,k)(α)]" = combo)
  colnames(out1) <- col_labs
  out2<-t(out1)
  return(out2)
}

## ========================================================================================
##  Examples
## ========================================================================================
## Example (1) : Thin Cantor Set
EndpointCalculator(n=2, p=1, q=3)
#      a_(n,k)(α) b_(n,k)(α) [a_(n,k)(α),b_(n,k)(α)]
# k=0 "0/9"      "9/9"      "[0/9,9/9]"            
# k=1 "0/9"      "1/9"      "[0/9,1/9]"            
# k=2 "2/9"      "3/9"      "[2/9,3/9]"            
# k=3 "6/9"      "7/9"      "[6/9,7/9]"            
# k=4 "8/9"      "9/9"      "[8/9,9/9]" 

## Example (2) : Thick Cantor Set
EndpointCalculator(n=2, p=1, q=4)
#      a_(n,k)(α) b_(n,k)(α) [a_(n,k)(α),b_(n,k)(α)]
# k=0 "0/16"     "16/16"    "[0/16,16/16]"         
# k=1 "0/16"     "2.5/16"   "[0/16,2.5/16]"        
# k=2 "3.5/16"   "6/16"     "[3.5/16,6/16]"        
# k=3 "10/16"    "12.5/16"  "[10/16,12.5/16]"      
# k=4 "13.5/16"  "16/16"    "[13.5/16,16/16]"  
\end{verbatim}

\subsection*{Appendix.C. Glossary}\label{Appendix.C}

This appendix includes notational conventions and key terms used throughout the survey for the convenience of readers outside descriptive set theory and fractal geometry. Specialists may skip this subsection without loss.\par 
The terms in the paper that a general mathematical reader (e.g., an analyst or topologist) may pause  fall into four clusters:\par 

\begin{enumerate}
    \item \textbf{Topological / set-theoretic descriptors} --- \textit{perfect set, nowhere dense, totally disconnected, derived set, $F_\sigma$, $G_\delta$, Borel hierarchy.}  
    
    \item \textbf{Measure-theoretic descriptors} --- \textit{Lebesgue measure, outer content, measure-zero vs. positive measure, self-similar measure, Cantor-Lebesgue function, Devil's staircase.}
    
    \item \textbf{Fractal / dynamical terminology} --- \textit{Hausdorff dimension, attractor, IFS, contraction coefficient, affine map, thickness.}
    
    \item \textbf{Continuum-theoretic terms} --- \textit{continuum, arc, dendrite, Peano continuum, embedding}  
\end{enumerate}
The concise mathematical definitions and descriptions are as in follows:\par 

\begin{itemize}
    \item \emph{Perfect set} --- a closed set in which every point is a
    limit point of the set; equivalently, a closed set with no isolated
    points~\cite{Munkres2000}.
    \item \emph{Nowhere dense set} --- a set whose closure has empty
    interior; informally, a set that is ``full of holes'' in every
    sub-interval.
    \item \emph{Totally disconnected set} --- a set whose only connected
    subsets are singletons; equivalently, between any two distinct points
    one can interpose a gap.
    \item \emph{$F_\sigma$ set} --- a countable union of closed sets.
    \item \emph{$G_\delta$ set} --- a countable intersection of open sets.
    \item \emph{Borel hierarchy} --- the transfinite stratification
    $\Sigma^0_n,\,\Pi^0_n$ of Borel sets by alternating countable unions
    and intersections, beginning with open and closed sets~\cite{Kechris1995}.
    \item \emph{Affine map on $[0,1]$} --- a function of the form
    $T(x) = rx + s$; a \emph{contraction} if $|r| < 1$, with
    $r = \mathrm{cont.\,coeff}(T)$ its \emph{contraction coefficient}.
    \item \emph{Iterated function system (IFS)} --- a finite family
    $\{T_k\}_{k=1}^{K}$ of contractions; its \emph{attractor} is the
    unique non-empty compact set $\Gamma$ satisfying
    $\Gamma = \bigcup_{k=1}^{K} T_k(\Gamma)$~\cite{Falconer2014}.
    \item \emph{Hausdorff dimension} --- the critical exponent
    $\dim_H(E)$ at which the $s$-dimensional Hausdorff measure of $E$
    transitions from $+\infty$ to $0$; for self-similar Cantor sets
    generated by $K$ maps with common ratio $r$, $\dim_H(E) = \log K /
    \log(1/r)$.
    \item \emph{Self-similar measure} --- a Borel probability measure
    $\mu$ on the attractor of an IFS $\{T_k\}$ satisfying
    $\mu = \sum_{k=1}^{K} w_k\,(\mu \circ T_k^{-1})$ for a probability
    weight vector $(w_1,\dots,w_K)$.
    \item \emph{Cantor--Lebesgue function (Devil's staircase)} --- the
    continuous, non-decreasing function $[0,1] \to [0,1]$ that is
    constant on every gap of the middle-third Cantor set and increases
    only on the set itself.
\end{itemize}


\end{document}